\def\submission{1}  

\ifdefined\submission
	\documentclass[preprint,review,12pt]{elsarticle}  
\fi
\ifdefined\final
	\documentclass[final,5p,times,twocolumn]{elsarticle}  
\fi

\usepackage{algorithmicx}
\usepackage{array}
\usepackage{amsmath,amssymb,amsthm}
\usepackage{booktabs}
\usepackage{graphicx,stfloats}
\usepackage{lineno}
\usepackage{multirow}
\usepackage[version=4, arrows=pgf]{mhchem}
\usepackage{url}
\usepackage{xspace}
\usepackage{hyperref}
\usepackage{hyphenat}

\usepackage{algpseudocode,algorithm}

\biboptions{sort&compress}
\graphicspath{{figs/}}
\modulolinenumbers[5]

\journal{Combustion and Flame}

\newcommand{\ie}{\textit{i.e.},\xspace}
\newcommand{\eg}{\textit{e.g.},\xspace}

\newcommand{\Ns}{N_{\mathrm{s}}}
\newcommand{\Nr}{N_{\mathrm{r}}}
\DeclareMathOperator*{\minimize}{\text{minimize}}

\DeclareMathOperator*{\st}{\text{subject to}}

\def\hyph{-\penalty0\hskip0pt\relax}

\begin{document}

\begin{frontmatter}
\title{A New Data-Driven Sparse-Learning Approach to Study Chemical Reaction Networks}

\author[fir]{Farshad Harirchi}
\author[sec]{Doohyun Kim}
\author[fir]{Omar Khalil}
\author[fir]{Sijia Liu}
\author[sec]{Paolo Elvati}
\author[fir]{Alfred Hero\corref{cor1}}
\ead{hero@eecs.umich.edu}
\author[sec,ter]{Angela Violi\corref{cor1}}
\ead{avioli@umich.edu}

\address[fir]{Department of Electrical Engineering and Computer Science, University of Michigan, Ann Arbor, MI 48109--2125, USA}
\address[sec]{Department of Mechanical Engineering, University of Michigan, Ann Arbor, MI 48109--2125, USA}
\address[ter]{Departments of Chemical Engineering, Biomedical Engineering, Macromolecular Science and Engineering, Biophysics Program, University of Michigan, Ann Arbor, MI 48109--2125, USA}
\cortext[cor1]{Corresponding authors:}

\begin{abstract}
Chemical kinetic mechanisms can be represented by sets of elementary reactions that are easily translated into mathematical terms using physicochemical relationships. The schematic representation of reactions captures the interactions between reacting species and products. 
Determining the minimal chemical interactions underlying the dynamic behavior of systems is a major task. 
In this paper, we introduce a novel approach for the identification of the influential reactions in chemical reaction networks for combustion applications, using a data-driven sparse-learning technique. 
The proposed approach identifies a set of influential reactions using species concentrations and reaction rates,
with minimal computational cost without requiring additional data or simulations.
The new approach is applied to analyze the combustion chemistry of \ce{H2 and C3H8} in a constant-volume homogeneous reactor. The influential reactions identified by the sparse-learning method are consistent with the current kinetics knowledge of chemical mechanisms. 
Additionally, we show that a reduced version of the parent mechanism can be generated as a combination of the influential reactions identified at different times and conditions and that for both \ce{H2 and C3H8} this reduced mechanism performs closely to the parent mechanism as a function of ignition delay over a wide range of conditions.
Our results demonstrate the potential of the sparse-learning approach as an effective and efficient tool for mechanism analysis and mechanism reduction. 
\end{abstract}


\end{frontmatter}

\section{Introduction}
\label{sec:introduction}
Combustion is a ubiquitous process with applications that vary extensively, ranging from heating and transportation to mass production of metallic and ceramic nanoparticles.
Combustion involves different processes that span a wide range of length and time scales, and are controlled by the delicate interplay among several chemical and physical phenomena, such as chemical kinetics, thermodynamics, fluid mechanics, and heat transfer. 
Over the past decades, chemical modeling has played a great role due to the continuing expansion in computational capabilities, development of diagnostic methods, quantum chemistry methods and kinetic theory. 
Within this context, the approach to developing detailed kinetic models (mechanisms) of fuel combustion involves compiling a set of elementary reactions whose rate parameters may be determined from individual rate measurements, reaction-rate theory, or a combination of both. 
These detailed mechanisms consist of a large number of chemical species and reactions: for example, the mechanisms of large hydrocarbons, as the ones used in surrogates for real transportation fuels, typically describe the interactions of $\sim 10^2$ -- 10\textsuperscript{3} species via 10\textsuperscript{3} -- 10\textsuperscript{4} reactions~\cite{law2007combustion}. 
The resulting networks of reactions show high nonlinearity and dimensionality. 
The relevance of the different species and reactions not only changes with time during combustion, but  also depends on the specific system and conditions (\eg temperature, pressure, equivalence ratio). 

The analysis of chemical kinetic mechanisms, particularly understanding the relevance of different species and reactions in various conditions, has several important applications.
While chemical mechanisms have been active research area for last several decades, refinements in implemented pathways and reaction rate parameters are still necessary. Thus, identification of most influential reactions in the network provides a critical information on which chemistry in the mechanism should be improved.
The analysis of chemical reaction networks can be leveraged to generate reduced chemical mechanisms that have similar combustion characteristics to the parent mechanism but present a significantly smaller number of species and reactions.
A direct application of this work is the use of reduced mechanisms into computational fluid dynamics simulations to describe realistic combustion systems like engines, overcoming the complexity and computational cost of large kinetic mechanisms that make full kinetic models unusable in most cases~\cite{tomlin_chapter_1997,turanyi_sensitivity_1990,lu2009toward}.
For example, computational simulations of internal combustion engines require the coupling of chemical models for the conversion of the fuel into combustion products with numerical treatments of the fluid dynamics of reacting flows~\cite{law2007combustion}.
Evaluating the influence of each species or reaction on the combustion process is a widely used strategy for developing reduced kinetic mechanisms~\cite{tomlin_chapter_1997,turanyi_sensitivity_1990,vajda_principal_1985,maas_implementation_1992,maas_simplifying_1992,lam_csp_1994,lam_using_1993,bhattacharjee_optimally-reduced_2003,lu_directed_2005,lu_strategies_2008,pepiotdesjardins_efficient_2008,niemeyer_skeletal_2010}.
Perturbation sensitivity analysis is to date the most commonly used method.
Of particular interest is the influence of kinetics parameters on combustion characteristics, such as ignition delay or heat release trends. 
The brute force approach to such sensitivity analysis is to evaluate the impact of each reaction on the combustion feature of interest by perturbing the rates one by one~\cite{chang_computational_1987,lifshitz1980oxidation,fridlyand2017role}. 
This brute force scheme is straightforward to implement but becomes overly computationally expensive for large mechanisms as a full simulation run is required for every change in the kinetic parameters.
More details on various sensitivity analysis methods and their applications can be found in~\cite{tomlin_chapter_1997,turanyi_sensitivity_1990}. 

An alternative approach to sensitivity analysis uses the rates of production of the species.
After a simulation, the reaction rate of each reaction in the original solution is analyzed to find reactions that contribute to production/consumption of some of the key species.
While this method can be effective for discovery of influential species, the method is not effective for general combustion system analysis across different thermodynamic conditions. 

In addition to these two common techniques, various other reaction-based methods can be found in the literature that aim at identifying unimportant reactions in the chemical network~\cite{vajda_principal_1985,maas_implementation_1992,maas_simplifying_1992,lam_csp_1994,lam_using_1993,bhattacharjee_optimally-reduced_2003}. 
These methods are often presented as mechanism reduction approaches.
However, in general, the computational cost of these techniques is high when applied to large scale  mechanisms, as many perturbations of the rate parameters or Jacobian evaluation of functions are required.

In this paper, we propose a novel approach based on data-driven sparse-learning methods to analyze reaction mechanisms and identify the most influential reactions. 
The data-driven approach simplifies complex nonlinear dynamical models to scalable linear models; the sparse-learning methodology removes weak interdependencies with the reaction network and creates a sparse network. 
Unlike other methods in the literature, the proposed method does not require any understanding of the underlying chemical process, instead it learns from data generated from the chemical process. 
The proposed method is applied to the analysis of two mechanisms describing the chemistry of \ce{H2} fuel and \ce{C3H8} fuel, respectively. 
We also demonstrate how this new method can be used to derive a reduced version of the parent chemical mechanism without the need for additional simulations, making it considerably more efficient compared to other equivalent methods. 
By finding a set of influential reactions using species concentrations and reaction rates, the proposed method is not biased toward any specific combustion property \cite{Harirchi2017Data}.


\section{Methodology}
\label{methodology}

\subsection{Model Framework}
The time evolution of a reactive system is commonly modeled using mass-action kinetic equations~\cite{feinberg87,gunawardena03,gillespie07,Chellaboina09,anderson2011continuous}, which, by applying Euler discretization method~\cite{ascher98}, can be written as
\begin{equation}
  \label{eqn:discrete_model}
  \vspace*{-0.2cm}
  \mathbf{X}_{t+1} = \mathbf{X}_t + \mathbf M \mathbf r_t \Delta_t + \pmb {\omega}_t, \; t = 1,2, \hdots , T-1,
\end{equation}
where $\mathbf{X}_t $ is an $\Ns \times 1$ vector of concentrations of all the $\Ns$ species at time $t$, $\Nr$ is the number of reactions whose stoichiometric coefficients are contained in the $\Ns \times \Nr$ matrix $\mathbf{M}$, $\Delta_t$ is the sampling time used for discretization, and $\pmb {\omega}_t$ denotes the discretization error and process noise.
Finally, $\mathbf{r}_{t}$ is the $\Nr \times 1$ vector of the rates of $\Nr$ reactions, which can be expressed as:
\begin{equation}
  \mathbf{r}_{t}(i) = k_i\prod_{j=1}^{\Ns}\mathbf{X}_t(j)^{\nu_{ij}}, 
\end{equation}
where $k_i$ is the rate constant for the $i^{th}$ reaction, $\mathbf{X}_t(j)$ denotes the $j^{th}$ component of $\mathbf{X}_t$, and $\nu_{ij}$ are the stoichiometric coefficients for the reactants of reaction $i$ and zero otherwise.

For each reaction, we introduce a selection (binary case) or weight (real-valued case) variable $w_{i,t}$, which describes the significance of the $i^{th}$ reaction at time $t$ (larger values correspond to more influential reactions), and the selection/weight vector $\mathbf w_t = [\mathrm w_{1,t}, \mathrm w_{2,t}, \ldots, \mathrm w_{\Nr, t}]^\intercal$, which represents the weights of all the reactions at time $t$. 

Considering the binary case, where $\mathbf w_t$ selects a subset of all reactions to describe the evolution of the system, an error is introduced in the concentration of the species that can be expressed for the $j^{th}$ species as
\begin{equation}
\mathcal{E}_{j,t}(\mathbf w_t) = |\mathbf X_{t+1}(j)-\mathbf X_t(j) - \mathbf M_{j} ( \mathbf w_t \odot \mathbf r_t ) \Delta_t|,
\end{equation}
where $\mathbf M_j$ denotes the $j^{th}$ row of matrix $\mathbf M$ and $\odot$ denotes the element-wise product. Based on these definitions, we can impose a constraint on the concentration error of each species at all times as
\begin{equation} 
  \label{eqn:err_tol}
  \mathcal{E}_{j,t}(\mathbf w_t) \leq \epsilon\mathcal{N}_{t}(j),
  ~ \forall t, \; \forall j \in \{1, \dotsc, \Ns \},
\end{equation}
\begin{equation}
  \label{eqn:norm}
  \mathcal{N}_t(j) = \lvert \mathbf M_j \rvert \mathbf r_t \Delta_t, \; j \in \{1, \dotsc, \Ns \}.
\end{equation}
where $\mathcal{N}_t$ is a normalization factor at time $t$, defined as the sum of absolute changes in all concentrations at time $t$, and $\epsilon$ is a tuning parameter that indicates the acceptable error tolerance of $\mathcal{N}_t$, \eg $\epsilon=0.05$ enforces a maximum of 5\% error.

The constraint introduced in Eq.~\ref{eqn:err_tol} is effective against the addition of noise in the concentration evolution but is not very effective in limiting constant drifts in the species concentration.
To correct for this issue, we added an additional constraint to limit the propagation of error over time.
For the change in concentration of the $j^{th}$ species in time horizon $\{t, t+1,\ldots,t+H-1\}$, we have that
\begin{equation}
  \label{eqn:err_prop}
  \lvert \mathbf X_{t+H-1}(j)- \mathbf{X}_{t}(j) - \sum_{k=t}^{t+H-2} \mathbf{M}_j (\mathbf{w}_{k}\odot \mathbf r_k) \Delta_k \rvert \leq \beta \epsilon \sum_{k=t}^{t+H-2} \mathcal{N}_{k}(j),
\end{equation}
where $H$ is the number of time samples $\{t, t+1, \dotsc\}$ in the time horizon $[t, t+H-1]$ and $\beta$ is a tuning parameter that limits the amount of concentration drift. For the choice of $H$, two points should be considered:
\begin{enumerate}
    \item \textit{Physical effect}: larger choice of $H$ enforces the aggregated error in the concentrations of all species that is caused by removing some of the reactions to remain in the acceptable range for a larger time horizon. On the other hand, the choice of smaller $H$, relaxes the propagated error to be in the same range , but for a smaller time horizon. Consequently, the larger $H$ results in a reduced mechanism with more reactions.   
    \item \textit{Computational cost}: larger $H$ results in less number of optimization problems of larger size, and smaller choice of $H$ creates more optimization problems of smaller size. The effect of $H$ on the number of variables in the optimization is illustrated in Table \ref{tab:complexity}.  
\end{enumerate}

\subsection{Sparse-Learning Reaction Selection (SLRS) method}
The task of finding the most influential reactions is a matter of determining the smallest subset of reactions that are active at any given time, such that the error in the concentrations induced by limiting the number of reactions remains in a user-specified tolerance range $(\epsilon, \beta)$ at all times.

To solve this problem, we formulated the following mathematical approach to data-driven sparse-learning reaction selection.
For each time batch of size $H$, we solve the following integer linear programming (ILP) problem:
\begin{align}
\tag{$P_{SLRS}$}
\displaystyle \{ \mathbf w_k^{\ast}\}_{k=t}^{t+H-1} = \minimize_{\{ \mathbf w_k\}_{k=t}^{t+H-1}} & \sum_{k=t}^{t+H-1}\sum_{i=1}^{\Nr} \mathbf w_{i,k} \label{eq:prob_ori}\\
  \st & \;\; \mathbf w_{i,k} \in \{0,1\}, ~i=1, 2, \dotsc, \Nr, \nonumber \\
      & \;\; \forall k \in \{t,t+1,\hdots ,t+H-1\}  \nonumber \\
      & \;\; \text{Eqs.~\ref{eqn:err_tol},\ref{eqn:err_prop} hold.}  \nonumber
\end{align}
where $\{\mathbf w_k\}_{k=t}^{t+H-1}$ are binary-valued optimization variables.
Note that solving problem~\eqref{eq:prob_ori} delivers the minimum number of reactions such that the error tolerance constraints on individual concentrations (Eq.~\ref{eqn:err_tol}) and on error propagation (Eq.~\ref{eqn:err_prop}) are satisfied.

To further reduce the complexity of our data-driven sparse-learning approach, we can use convex relaxation methods~\cite{NoceWrig06} by replacing the binary variable constraint with its convex \textcolor{black}{hull}, \ie $0 \leq \mathbf w_{i,t} \leq 1$ for all $i$.
The relaxed problem has the following form:
\begin{align}
\tag{$P_{RSLRS}$}
\displaystyle \{ \mathbf w_k^{\ast}\}_{k=t}^{t+H-1}= \minimize_{\{ \mathbf w_k\}_{k=t}^{t+H-1}} & \sum_{k=t}^{t+H-1}\sum_{i=1}^{\Nr} \mathbf w_{i,k} \label{eq:prob_rel}\\
 \st &  \;\; 0 \leq \mathbf w_{i,k} \leq 1, ~i=1,2,...,N_r,   \nonumber \\
&  \;\; \forall k \in \{t,t+1,\hdots ,t+H-1\}  \nonumber \\
&  \;\; \text{Eqs.~\ref{eqn:err_tol},\ref{eqn:err_prop} hold.} \nonumber
\end{align}

This relaxation of \eqref{eq:prob_ori} yields a linear programming problem, which can be solved efficiently in polynomial time, and therefore makes our approach promising for large-scale chemical reaction networks at the cost of losing optimality guarantees.
Moreover, the solution of~\eqref{eq:prob_rel} $\mathbf{w}_t^{\ast}$ is a real-valued vector, which is more informative but not immediately suitable for the selection of a subset of reactions.
In order to project this solution onto a binary vector, we set a threshold as follows
\begin{equation}
\begin{cases}
  \tilde{\mathbf{w}}_{i,k} = 1 \qquad \text{if} \; \mathbf{w}^{\ast}_{i,k} > \alpha \\
  \tilde{\mathbf{w}}_{i,k} = 0 \qquad \text{if} \; \mathbf{w}^{\ast}_{i,k} \leq \alpha \\
\end{cases},
\label{eqn:cutoff}
\end{equation}
where $\alpha$ is a threshold value and $\tilde{\mathbf{w}}$ is in the form of a selection vector similar to the solution of problem~\eqref{eq:prob_ori}.

The next step is to find the reduced mechanism not just for the time horizon $[t,t+H-1]$, but for the entire sampled interval $[0,t_{f}]$, where $t_f$ is the last time sample.
In order to do this, we solve problem~\eqref{eq:prob_rel} for time intervals $[0,H-1], [H,2H-1], \ldots,[KH,t_f]$, where $KH$ is the largest integer multiple of $H$ that is smaller than or equal to $t_f$, as described in Algorithm~\ref{alg:1}.
The output of Algorithm \ref{alg:1} is a matrix $\mathbf W(\theta)$, with columns indicating the selected reactions at each time instance under specific initial conditions, namely temperature $T$, equivalence ratio $\phi$, and pressure $P$ (\ie $\theta = [T,\phi,P]$).

\begin{algorithm}[!htb]
	\caption{Calculating $\mathbf{w}(\theta)$}
  \label{alg:1}
	Input: $\{\mathbf{X}_t\}_{t=0}^{t_f}$, $\;\{\mathbf{r}_t\}_{t=0}^{t_f}$, $\;\{\Delta_t\}_{t=1}^{t_f}$, $\;\epsilon$, $\;\beta$, $\;\theta = [\phi,T,P]$, H, $\mathbf{M}$\\
    Initialize:  $t=0$, $\mathbf{w}_{-1} = \mathbf 0$, $\mathbf W (\pmb{\theta}) \in \mathbb{R}^{N_r \times N_t}$.
	\begin{enumerate}
        \item While ($t \leq t_f$)
        \begin{itemize}
        	\item $W = \min\{ H, t_f-t+1 \}$.
            \item set the time horizon to $[t, t+W-1]$.
            \item solve problem~\eqref{eq:prob_rel} to obtain $\mathbf w^{\ast}_{t},t\in [t, t+W-1]$.
            \item Assign $\mathbf w^{\ast}_{t},t\in [t, t+W-1]$ values to columns [t,t+W-1] of $\mathbf W (\pmb{\theta})$.
			\item $t+W \rightarrow t$.
		\end{itemize}
        End while
        \item Return $\mathbf W (\pmb{\theta})$ or $\mathbf w^{\ast}_{t},t\in [0,t_f]$.
	\end{enumerate}
\end{algorithm}

\subsection{Complexity Analysis and Computational Cost}
\label{subsec:complexityAnalysisAndComputationalCost}
The original problem~\eqref{eq:prob_ori} is a standard ILP problem, and thus can be solved using state-of-the-art ILP solvers such as Gurobi~\cite{gurobi} and CPLEX~\cite{cplex}.
Even though these solvers can solve problems with large numbers of integer variables relatively fast by employing branch and bounding algorithms, the worst-case complexity of ILP is exponential in the number of integer variables.
On the other hand, the relaxed problem~\eqref{eq:prob_rel} is a linear programming (LP) problem, which can be solved in polynomial time.
Table~\ref{tab:complexity} lists the factors (\ie variables and constraints) that affect the performance of the two approaches. 
\begin{table*}[hbt]
  \centering
  \caption{Complexity analysis for each optimization problem with time horizon of size $H$}
  \label{tab:complexity}
	\begin{tabular}{lcccc}
  \toprule
	 & Real      & Integer   & Linear     & Integral  \\
	 & variables & variables & constraints & constraints \\
  \midrule
	\ref{eq:prob_ori} & 0       & $N_r H$ & $N_rH$            & $2N_s(H+1)$ \\
	\ref{eq:prob_rel} & $N_r H$ & 0       & $2N_s(H+1)+3N_rH$ & 0 \\
  \bottomrule
  \end{tabular}
\end{table*}

The performance of the integer programming solution and the relaxed linear programming solution are compared in Tab.~\ref{tab:runtime}. 
These run-time values are calculated by taking the average of the times it takes to run Algorithm~\ref{alg:1} (with the corresponding formulation) for different values of $\theta$. 
The results show that solving the relaxed formulation \eqref{eq:prob_rel} significantly reduces the computational cost to less than $\frac{1}{5}$th of the integer formulation.
It should be noted that both methodologies can easily be accelerated via parallelization, as each initial condition $\theta$ can be analyzed independently.
Since the original problem \eqref{eq:prob_ori} is not computationally scalable, the rest of this paper focuses on the relaxed problem.

\begin{table}[hbt]
  \centering
  \caption{Comparison of computational performances of \ref{eq:prob_ori} and \ref{eq:prob_rel}.
  Average and total analysis time based on 216 Propane mechanism simulations with different initial conditions.}
  \label{tab:runtime}
  \begin{tabular}{lrr}
  \toprule
	 & Average (core s) & Total (core hr)\\
  \midrule
  \ref{eq:prob_ori} & 2273 & 136.38 \\
  \ref{eq:prob_rel} & 394  & 23.64 \\
  \bottomrule
  \end{tabular}
\end{table}

The total analysis time for 48 \ce{H2}/air system simulations using sparse learning reaction selection approach is 2.99 seconds. All the calculations are performed on a Linux based machine with a 2.1 GHz processor and 8 GB of RAM.

\section{Analysis of \texorpdfstring{\ce{H2}}{H2}/air system}
As the first system of investigation, we used the proposed relaxed SLRS approach to study the combustion of \ce{H2} in air.
Hydrogen combustion is chosen because it yields one of the smallest reaction networks in combustion (8 species and 62 reactions counting both forward and reverse processes~\cite{hong_improved_2011}) and is a well-studied system.

\begin{figure}[htb]
\centering
\includegraphics[width=13cm,trim={0 1cm 0 0},clip]{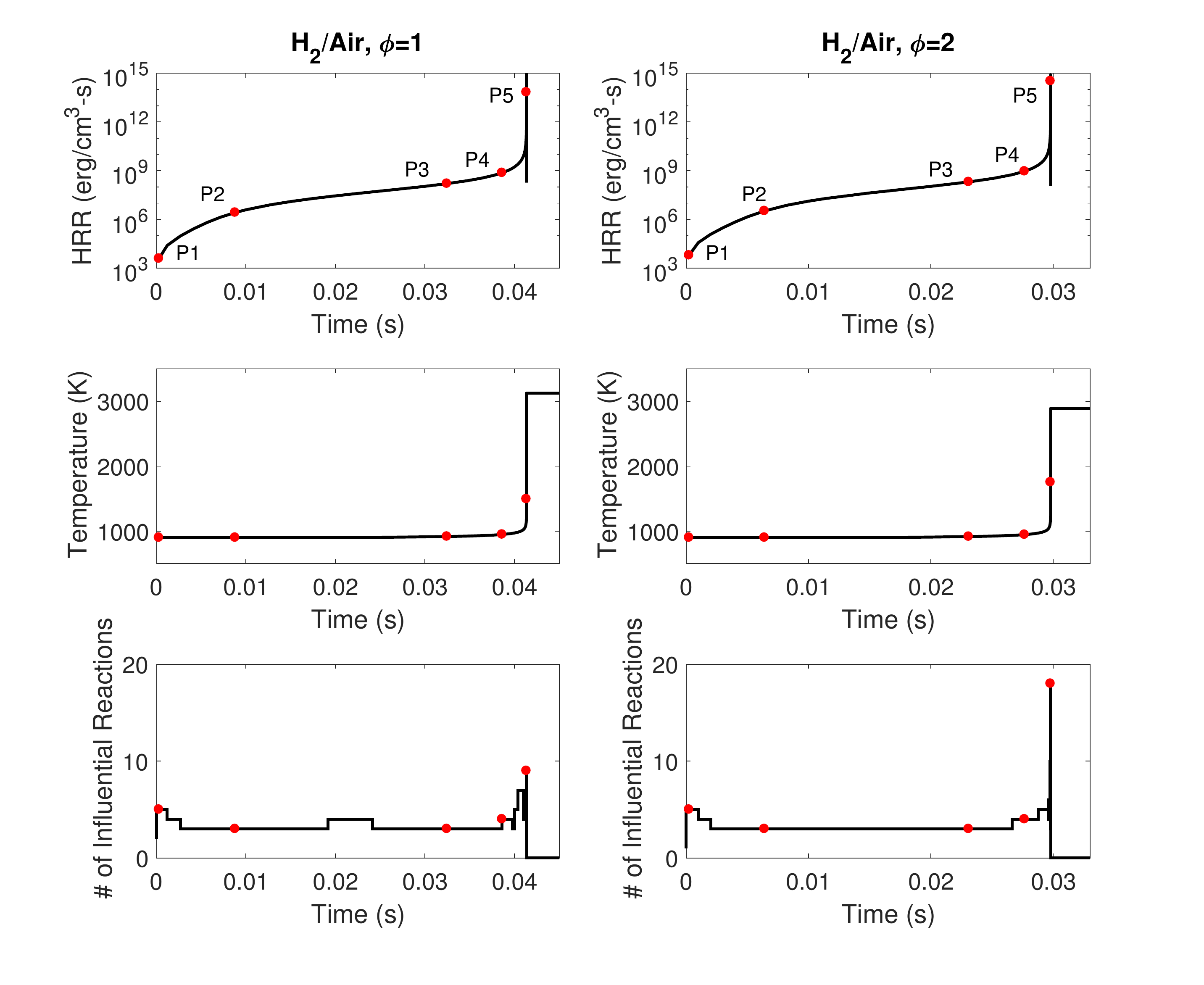}
\caption{Calculated heat release rates (HRR), temperature profiles, and the number of influential reactions for \ce{H2}/air mixture with the initial condition of 20~atm and 900~K. Highlighted are five points (P1-P5) during the time evolution. }
\label{fig:h2_phi1_2}
\end{figure}

Constant-volume, homogeneous reactor simulations, performed with Ansys Chemkin~\cite{noauthor_chemkin_2016} and a \ce{H2} oxidation mechanism by Hong et al.~\cite{hong_improved_2011}, were used to create the reference chemical reaction network for two different initial conditions ($P=20$~atm, $T=900$~K, $\phi=1,2$).
For each condition, we identified the influential reactions using SLRS parameters of $\epsilon=0.21$, $\beta=3$, and $\alpha=0$, and the results of five sampling regions are shown in Fig.~\ref{fig:h2_phi1_2}.

\begin{table}[hbt]
	\centering
	\caption{Influential reactions at selected time samples (see Fig.~\ref{fig:h2_phi1_2}) for \ce{H2}/air mixture with the initial condition of 20~atm and 900~K.
    Influential reactions ($w_{i,k} > 0$) are the same for $\phi=1$ and $\phi=2$ except at P5, for which, in fuel-rich conditions, additional 9 reactions are selected.}
    \label{tab:H2_phi1_2}
	\tiny	
    \begin{tabular}{cl@{}}
    \toprule
    Condition & Influential Reactions \\
    \midrule
 	P1           & \ce{H + O2 ({} + M) -> HO2 ({} + M)}\\  
    (all $\phi$) & \ce{H + O2 ({} + O2) -> HO2 ({} + O2)} \\
                 & \ce{H2 + HO2 -> H + H2O2} \\
                 & \ce{H2 + O2 -> H + HO2} \\
                 & \ce{H2 + OH -> H + H2O} \\ 

	\addlinespace[0.75em]
    P2           & \ce{H + O2 ({} + M) -> HO2 ({} + M)} \\  
    (all $\phi$) & \ce{H2 + HO2 -> H + H2O2} \\
                 & \ce{H2 + OH -> H + H2O} \\ 

	\addlinespace[0.75em]
    P3           & \ce{H + O2 ({} + M) -> HO2 ({} + M)} \\ 
    (all $\phi$) & \ce{H2 + OH -> H + H2O} \\ 
                 & \ce{HO2 + HO2 -> O2 + H2O2} \\

	\addlinespace[0.75em]
    P4           & \ce{H + O2 ({} + M) -> HO2 ({} + M)} \\
    (all $\phi$) & \ce{H2 + OH -> H + H2O} \\  
    			 & \ce{H2O2 + H -> H2O + OH} \\
                 & \ce{HO2 + HO2 -> O2 + H2O2} \\

	\addlinespace[0.75em]
 	P5 			& \ce{H + HO2 -> OH + OH} \\
   (all $\phi$) & \ce{H + HO2 -> H2 + O2} \\
    			& \ce{H + O2 -> O + OH} \\
                & \ce{H2 + O -> H + OH} (DUP) \\
                & \ce{H2 + OH -> H + H2O} \\ 
                & \ce{H2O2 ({} + M) -> 2OH ({} + M)} \\
                & \ce{H2O2 + H -> H2O + OH} \\
                & \ce{OH + HO2 -> H2O + O2} \\

	\addlinespace[0.75em]
    P5        	  & \ce{2H ({} + M) -> H2 ({} + M)} \\
  (only $\phi=2$) & \ce{2H + H2  -> 2H2} \\
                  & \ce{H + HO2  -> H2O + O} \\
                  & \ce{H + OH ({} + M) -> H2O ({} + M)} \\  
                  & \ce{H2O + O  -> OH + OH} \\
                  & \ce{H2O2 + H  -> HO2 + H2 } \\
                  & \ce{O + H ({} + M) -> OH ({} + M)} \\ 
                  & \ce{O + HO2  -> OH + O2 } \\
                  & \ce{OH + H + H2O  -> 2H2O} \\
    \bottomrule
	\end{tabular}
\end{table}

The analysis of the results, listed in Table~\ref{tab:H2_phi1_2}, shows that the set of influential reactions prior to the ignition time is identical for both stoichiometric and fuel-rich cases.
Initially at P1, the algorithm identifies reactions that consume the fuel as well as three-body reactions generating \ce{HO2} as influential.
As the ignition process progresses (from P2 to P3), reactions associated with \ce{HO2} production/consumption become prominent, with a shift in the \ce{H2O2} production pathway  from \ce{H2 + HO2} to recombination of two \ce{HO2} radicals, likely due to increased concentration of \ce{HO2} radicals.
As the system is approaching ignition (P4), the relaxed model introduces an additional influential reaction pathway that leads to the formation of \ce{OH}.

Once ignition is reached (P5), we observe a difference between the two systems in influential reactions. Interestingly, influential reactions from the stoichiometric case are a subset of the ones from the fuel-rich case.  
In both cases, the algorithm selects chain branching reactions (\eg \ce{H +O2 -> O + OH, \; H2 + O -> H + OH}) that are well-known factors for onset of high\hyph temperature ignition as well as reactions that are associated with \ce{HO2} and water formation.
In addition to these reactions, under fuel-rich conditions, nine more reactions are identified including three \ce{H2} formation reactions from \ce{H} radical, which describe the dynamic equilibrium between \ce{H} radical and remaining \ce{H2}.

The time evolution of the influential reactions identified by the proposed algorithm described above agrees well with the general understandings of the ignition of hydrogen as discussed in the literature~\cite{law2010combustion,hong_improved_2011}. 
For example, when the temperature is lower than 900~K (P1 to P4), the proposed algorithm correctly identifies the preferential pathway to produce \ce{HO2} from \ce{H + O2} instead of the chain branching pathway that generates \ce{O + OH}, which is important only at higher temperatures.
Moreover, the proposed algorithm correctly captures the dominance of fuel consumption reactions during the very early phase, as well as the reactions that lead to ignition at the ignition time (\ce{H2O2} dissociation, chain propagation/branching reactions for \ce{HO2}).

As the proposed sparse-learning approach is designed to identify the influential reactions of a combustion process, it can be leveraged to generate a reduced mechanism, \ie a coarser reaction network of the analyzed mechanism.
We emphasize that, while the reduction of the \ce{H2} mechanism is not intended for a real-world application, our results highlight the algorithm's capability to generate reduced mechanisms.

To this end, we generated and analyzed 48 homogeneous reactor constant-volume simulations of \ce{H2}/air combustion with initial conditions between 5--20~atm, 800--1100~K, and equivalence ratios between 0.5--2.
From these simulations, we selected the union of the influential reactions identified for each initial condition at each time, obtaining a set of 31 reactions.
The ignition delay times in constant-volume simulations for the reduced mechanism were compared against those obtained using the full mechanism.

\begin{figure}[ht]
\centering
\includegraphics[width=13cm,trim={0 0.5cm 0 0},clip]{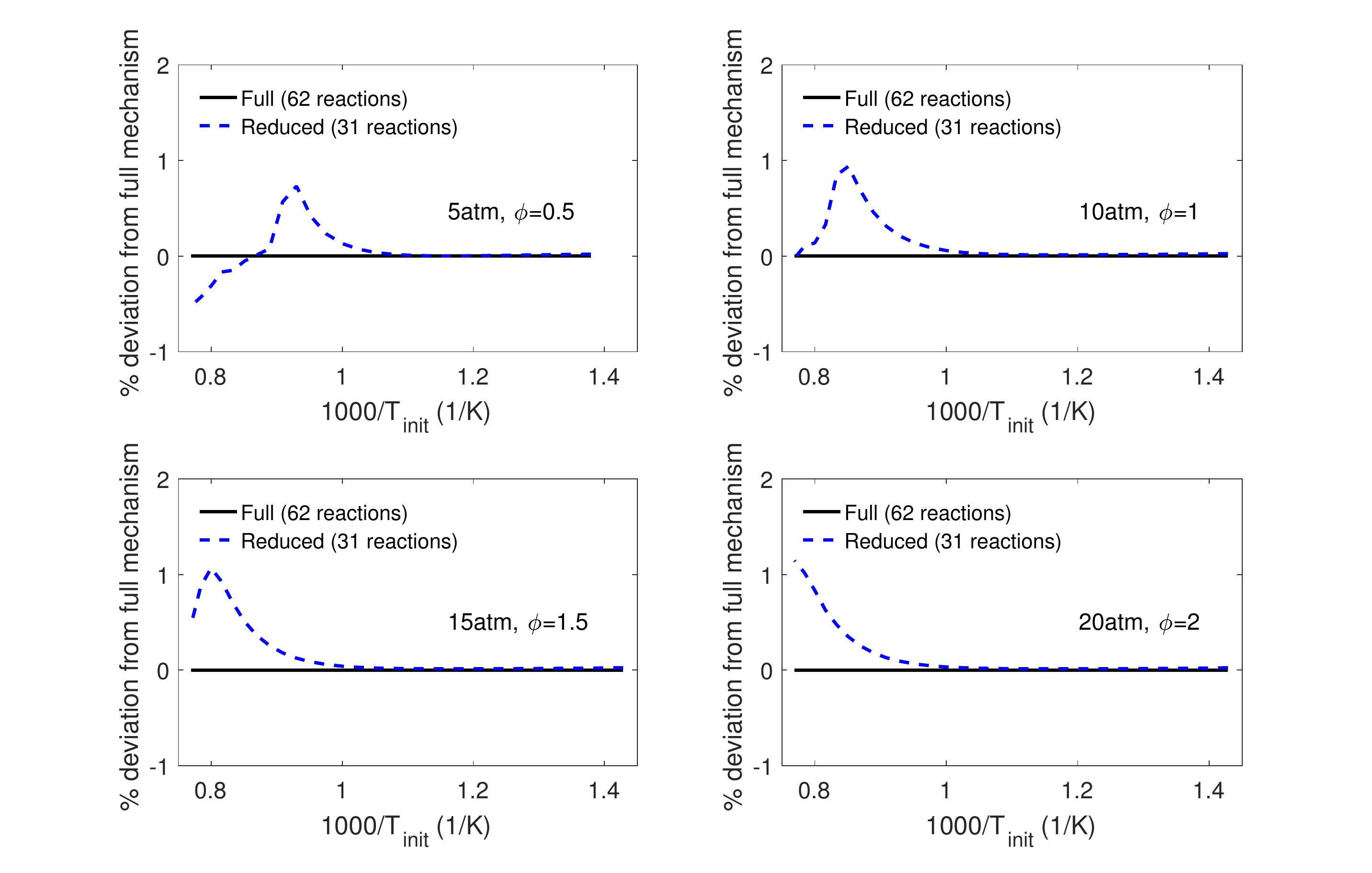}
\caption{Ignition delay time deviation from the full mechanism for \ce{H2}/air mixture. }
\label{fig:h2_ID}
\end{figure}

Figure~\ref{fig:h2_ID} shows the percentage deviation in ignition delay time from the full \ce{H2} mechanism for selected conditions.
The results which are representative of all the tested conditions, show that the differences in ignition delay time are well below 2\% for the reduced mechanism over a wide range of conditions, despite the significantly smaller number of reactions compared to the full mechanism (31 vs 62 reactions).
The maximum deviation occurs when the initial temperature is above 1000~K, while low\hyph temperature ignition delay times are nearly identical to the full mechanism. 
Overall, the excellent reproducibility of the ignition delay time (less than 2\% over all the conditions) indicates that the proposed method is a very effective approach to build a reduced mechanism that performs well in a wide range of thermodynamic conditions.
Moreover, it demonstrates that the ignition delay time, which is one of the most important combustion characteristics in practical energy conversion devices, can be preserved with the proposed method even if it is not explicitly taken into account directly during the reduction process.

\section{Analysis of Propane Combustion}

As a second application of the proposed sparse-learning approach, we analyzed the network of reactions describing the chemistry of propane. We used the full mechanism by Petersen et al.~\cite{petersen2007methane}, which includes 117 species and 1270 reactions.
Similar to the hydrogen case described in the previous section, we used a 0-D reactor for the simulation of a stoichiometric propane/air mixture at 20~atm and 700~K.
This system was chosen to identify the reactions that are responsible for the two-stage ignition behavior, a distinctive low\hyph temperature ignition characteristics of alkanes.

\begin{figure}[htb]
\centering
\includegraphics[width=6.5cm]{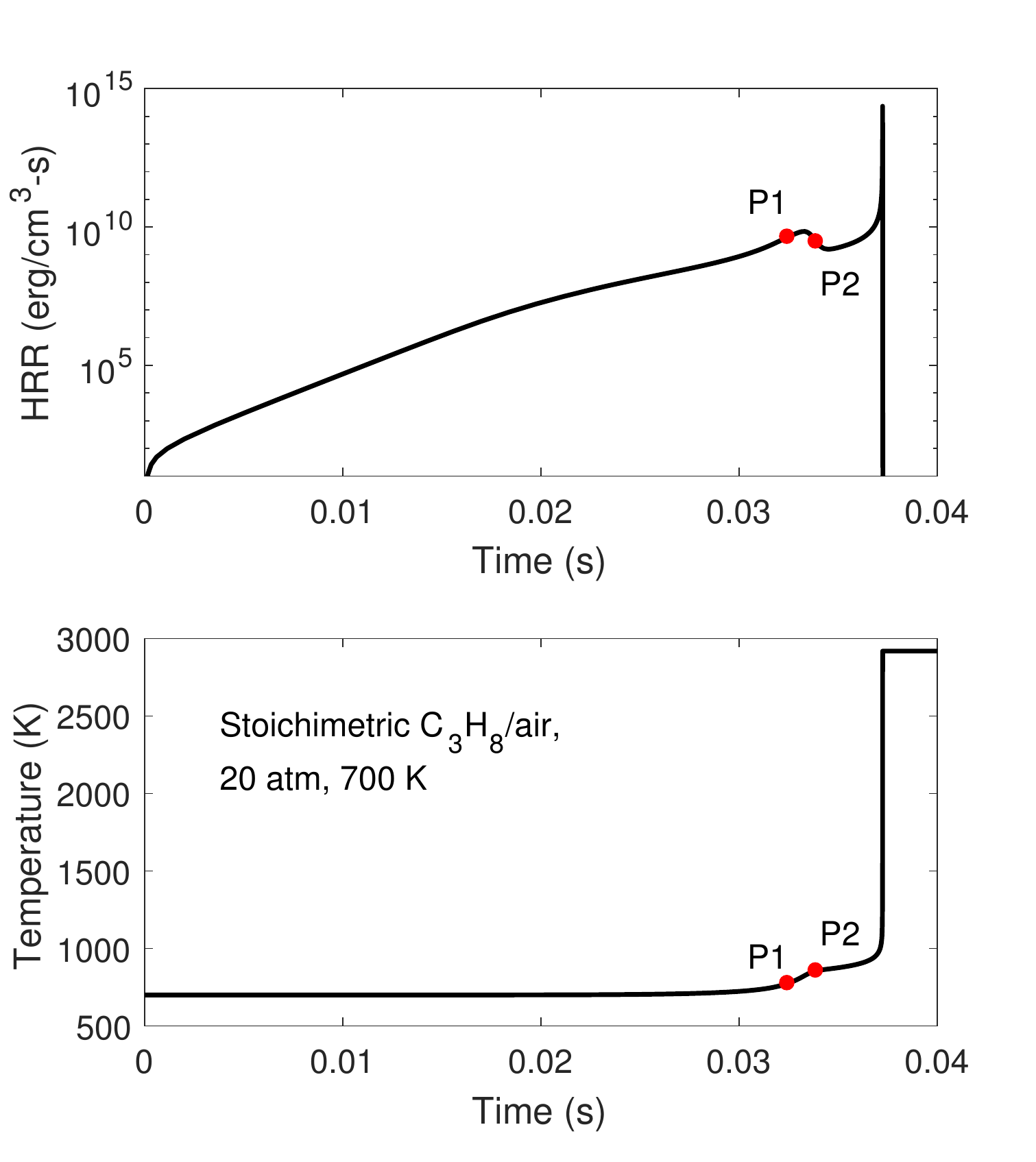}
\caption{Calculated heat release rate (HRR) and temperature of stoichiometric propane/air mixture with the initial condition of 20~atm and 700~K.}
\label{fig:propane_700K_phi1}
\end{figure}

\begin{figure}[htb]
\centering
\includegraphics[width=7.5cm]{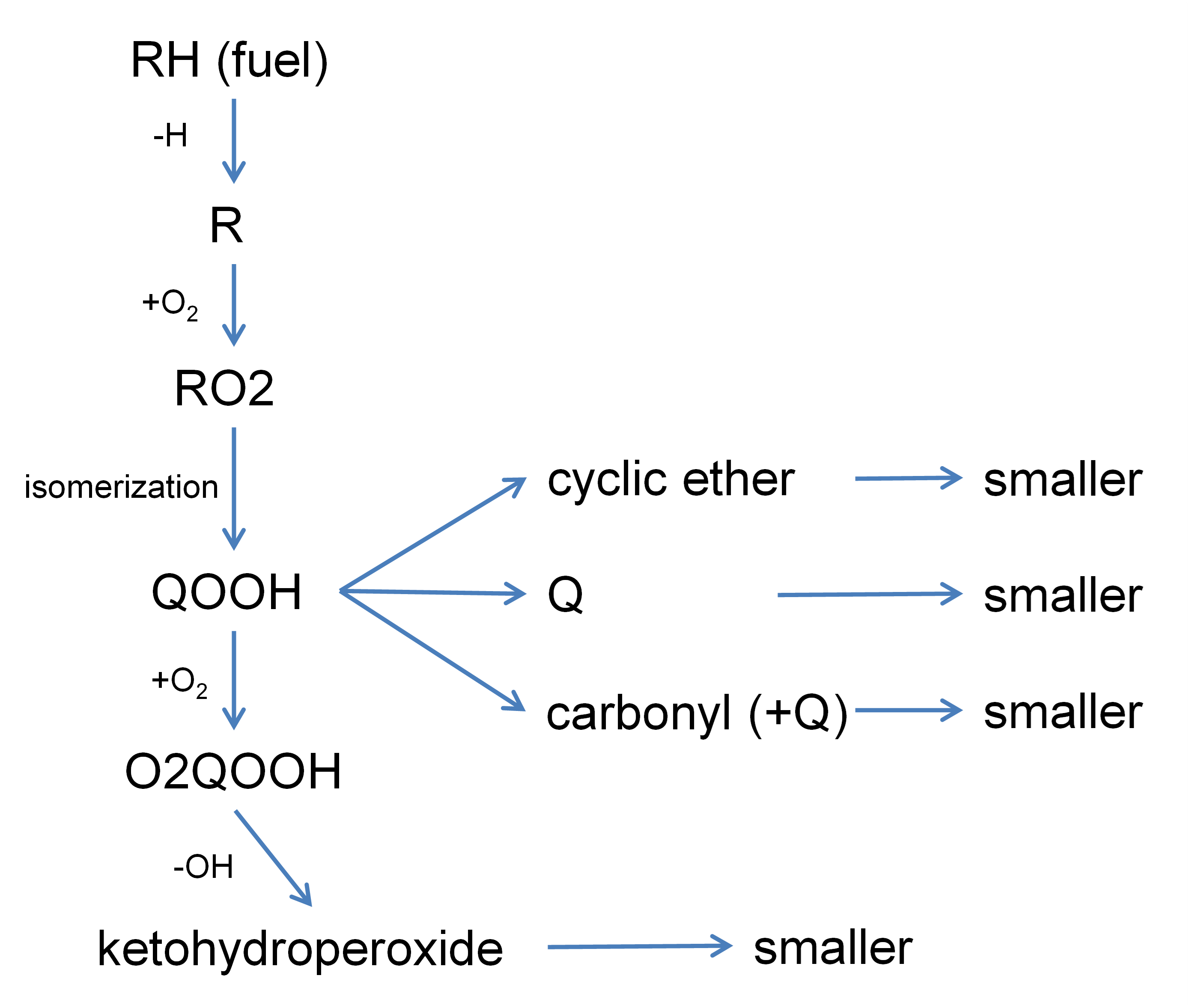}
\caption{Schematic of major low\hyph temperature oxidation pathways for alkanes. 
}
\label{fig:alkane_path}
\end{figure}

Fig.~\ref{fig:propane_700K_phi1} reports the calculated heat release rate and temperature profile for this system. The two-stage behavior is highlighted by the peak of heat release rate  and subsequent decrease (0.03 -- 0.035~s) before the ignition occurs ($\sim$~0.0372~s). This trend is generally associated with the competition between different low\hyph temperature pathways (see Fig.~\ref{fig:alkane_path} for a scheme of the major low\hyph  temperature pathways in alkane mechanisms~\cite{curran1998comprehensive,curran2002comprehensive,westbrook2011detailed,sarathy2011comprehensive}).

At low temperatures, once olefin peroxides (QOOH) are initially formed, the pathway to ketohydroperoxides via second \ce{O2} addition (\ce{QOOH ->} O2QOOH \ce{->} ketohydroperoxide) is favored.
However, as temperature increases during the ignition process, other pathways (\ie \ce{QOOH -> Q or carbonyl or cylic ether}) become increasingly relevant and dominate the low\hyph temperature reactivity.
As olefins (Q) and carbonyl radicals are more stable than ketohydroperoxides in such conditions, the charge reactivity and heat release rate slightly decrease as shown above.

To test if this behavior is captured by the proposed algorithm, we analyzed the system ($\epsilon=0.5$, $\beta=3$, $\alpha=0$) at two different times, namely, before (P1) and after (P2) the peak heat release rate point as shown in Fig.~\ref{fig:propane_700K_phi1}.
Out of 1270 reactions in the mechanism, the algorithm identified 68 influential reactions at P1 and 49 influential reactions at P2, with most of the differences between the two sets occurring in the low\hyph temperature chemistry reactions, which are listed in Table~\ref{tab:propane_phi1}.
The algorithm correctly identifies reactions in both groups of pathways as influential at P1, while at P2 the pathway to ketohydroperoxide is omitted and only the reactions that lead to Q or carbonyl production (less reactive pathway) are selected as influential.

\begin{table}[htb]
	\centering
	\caption{Low\hyph temperature chemistry subset of the influential reactions at selected times for the combustion of stoichiometric propane/air mixture with the initial condition of 20~atm and 700~K.
   P1 and P2 correspond to the instances shown in Fig.~\ref{fig:propane_700K_phi1}.
   Species names used in the mechanism~\cite{petersen2007methane} are shown.
   Q stands for olefin, and KET for ketohydroperoxide.
   Reaction classes  are labeled as in~\cite{sarathy2011comprehensive}.}
    \label{tab:propane_phi1}
  \small
  \begin{tabular}{@{}clc@{}}
  \toprule
  Condition & Influential reactions & Class\\
  \midrule
  P1 & C3H6OOH1-2 \ce{->} C3H6 + HO2       & 24 \\
    & C3H6OOH2-1 \ce{->} C3H6 + HO2        & 24 \\
    & C3H6OOH1-3 \ce{->} C3H6O1-3 + OH     & 25 \\
    & C3H6OOH2-1 + O2 \ce{->} C3H6OOH2-1O2 & 26 \\
    & C3H6OOH1-3 + O2 \ce{->} C3H6OOH1-3O2 & 26 \\
    & C3H6OOH1-2 + O2 \ce{->} C3H6OOH1-2O2 & 26 \\
    & C3H6OOH1-3O2 \ce{->} C3KET13 + OH    & 27 \\
    & C3H6OOH2-1O2 \ce{->} C3KET21 + OH    & 27 \\
    & C3H6OOH1-2O2 \ce{->} C3KET12 + OH    & 27 \\
    & C3KET13 \ce{->} CH2O + CH2CHO + OH   & 28 \\
    & C3KET21 \ce{->} CH2O + CH3CO + OH    & 28 \\
    & C3KET12 \ce{->} CH3CHO + HCO + OH    & 28 \\
  \addlinespace[0.75em]
  P2 & C3H6OOH2-1 \ce{->} C3H6 + HO2       & 24 \\
     & C3H6OOH1-2 \ce{->} C3H6 + HO2       & 24 \\
     & C3H6OOH1-3 \ce{->} C3H6O1-3 + OH    & 25\\
  \bottomrule
  \end{tabular}
\end{table}

Since our approach was able to capture the low\hyph temperature behavior of propane combustion, we analyzed the differences between low ($\phi=1$, 700~K, 20~atm, simulated above) and high\hyph temperature ($\phi=1$, 1500~K, 20~atm) combustion using the same parameters as before ($\epsilon=0.5$, $\beta=3$).
To simplify the comparison, we added an additional constraint to the optimization problem to accumulate the selection of influential reactions over the duration of each simulation.

\begin{figure*}[ht]
\centering
\includegraphics[width=14cm,trim={0 0.7cm 0 0},clip]{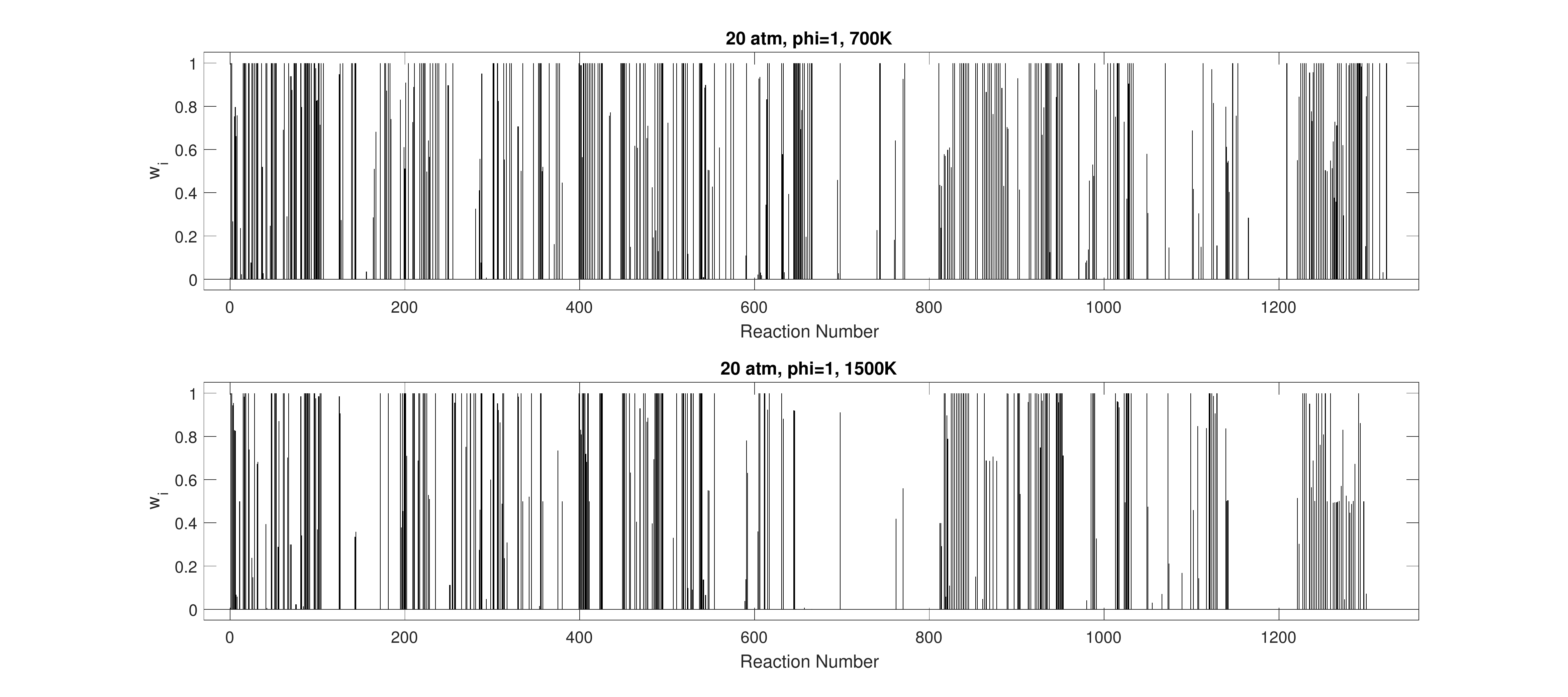}
\caption{Relevance (measured by the weight $w_i$) of the reactions from propane ignition simulations.
}
\label{fig:propane_700_1500K}
\end{figure*}

\begin{figure}[htb]
\centering
\includegraphics[width=12cm,trim={0 0.7cm 0 0},clip]{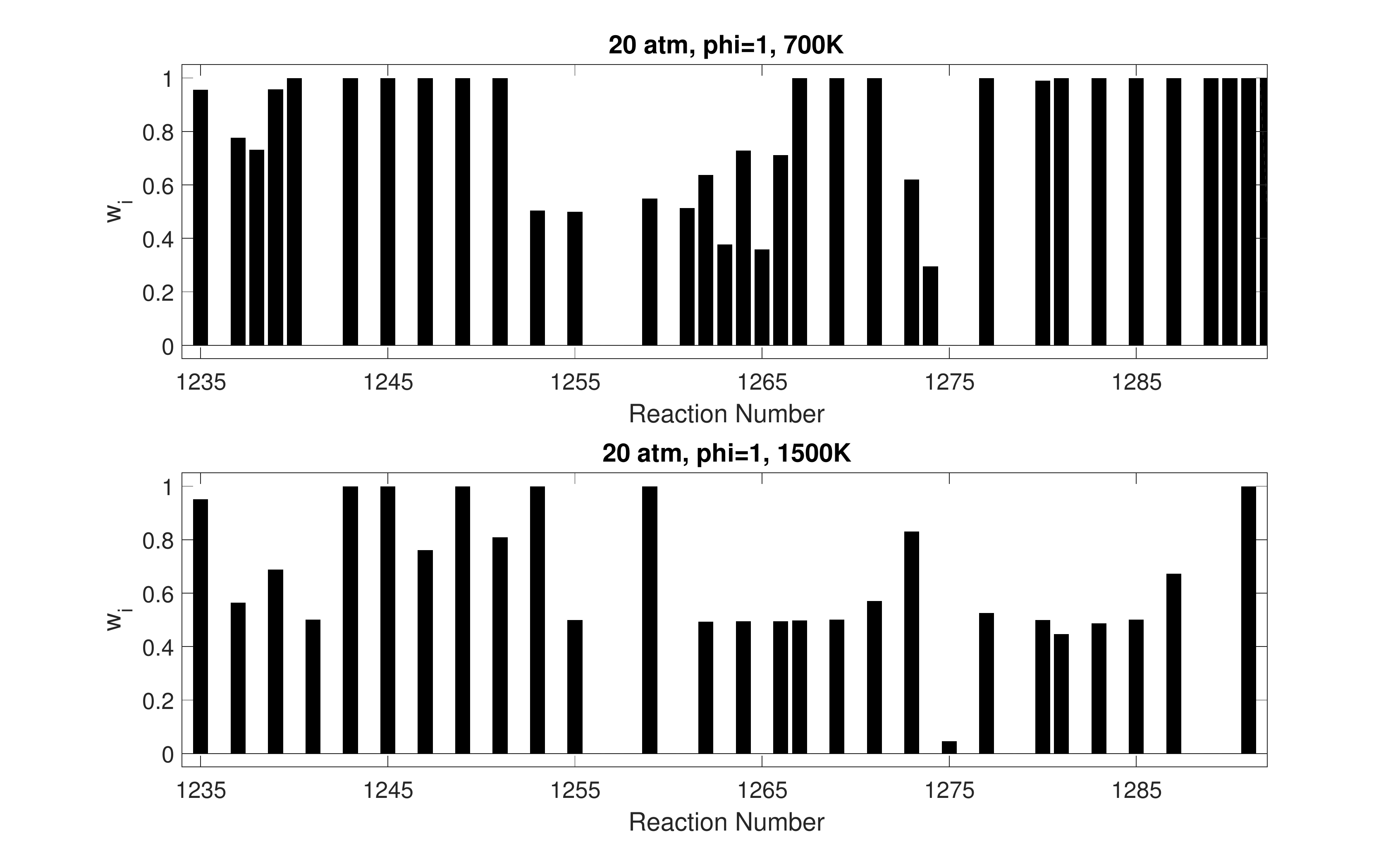}
\caption{Magnification of the importance (measured by the weight $w_i$) of the propane ignition reactions shown in Fig.~\ref{fig:propane_700_1500K}.
}
\label{fig:propane_700_1500K_low_temp}
\end{figure}

Figure~\ref{fig:propane_700_1500K} shows the relevance of each reaction in both conditions, as increasing values of $w_i$ correspond to more influential reactions.
The reactions are labeled by using an incremental number obtained by listing forward and reverse reactions in the same order they are presented in the original mechanism. 
For example, the first reaction in the parent mechanism is \ce{H + O2 <=> O + OH}, which becomes reaction 1 (\ce{H + O2 -> O + OH}) and 2 (\ce{O + OH -> H + O2}) in the abscissa of Figs~\ref{fig:propane_700_1500K} and Figs~\ref{fig:propane_700_1500K_low_temp}.
From the comparison of the two plots, we find that the algorithm identifies more influential reactions ($\alpha=0$) for the low\hyph temperature system (420) than for the high\hyph temperature one (350).  
This trend is expected, as the former system needs to include reactions from the low\hyph temperature regime in addition to the high\hyph temperature reactive pathways, which is in agreement with other works presented in the literature~\cite{niemeyer_skeletal_2010, ranzi2012hierarchical}.
The detailed comparison between the two sets of reactions can also provide useful insights on the low\hyph temperature chemistry.
For example, most of the reactions between 647--664, which describe hydrogen abstractions from \ce{C2H3CHO} by various radicals, are identified as influential only for the 700~K case.
The analysis of the mechanism indicates that \ce{C2H3CHO} is predominantly associated with the low\hyph temperature chemistry including \ce{RO2} and O2QOOH elimination pathways, \eg \ce{QOOH -> Q ->} smaller.
Similarly, reactions between 1301--1326, which represent decomposition pathway of C3 carbonyl compounds (C3H6O1-2, C3H6O1-3), are irrelevant for the high\hyph temperature system while some of them are influential for the 700~K case.
As these C3 carbonyls are mostly created by QOOH chemistry, pathways that are well-known to be associated with low\hyph temperature conditions, we find again that the algorithm selection matches our understanding of this kinetic mechanism.

To create a list of influential reactions a cutoff ($\alpha$) for $w_i$ (see Eq.~\ref{eqn:cutoff}) is required.
However, directly comparing the values of $w_i$ in their multireaction context is more informative. 
For example, if we analyze in detail the reactions between 1235--1288 (shown in  Fig.~\ref{fig:propane_700_1500K_low_temp}) that represent several low\hyph temperature reactions (\eg RO2 \ce{->} QOOH, QOOH destruction pathways, \ce{QOOH + O2 -> O2QOOH -> KET ->} smaller, O2QOOH \ce{->} smaller), we can see that while almost all the reactions have $w_i>0$, their relevance ($w_i$) changes in the two systems. 
In this way, it is straightforward to identify the reactions that are relevant in both regimes (\eg 1245) or reactions that are prevalent in a specific regime (\eg 1259 or 1267), particularly for complex and larger reaction networks, for which the assessments of relative importance among different pathways for specific combustion behavior is very challenging.

\begin{figure}[ht]
\centering
\includegraphics[width=13cm,trim={0 0.7cm 0 0},clip]{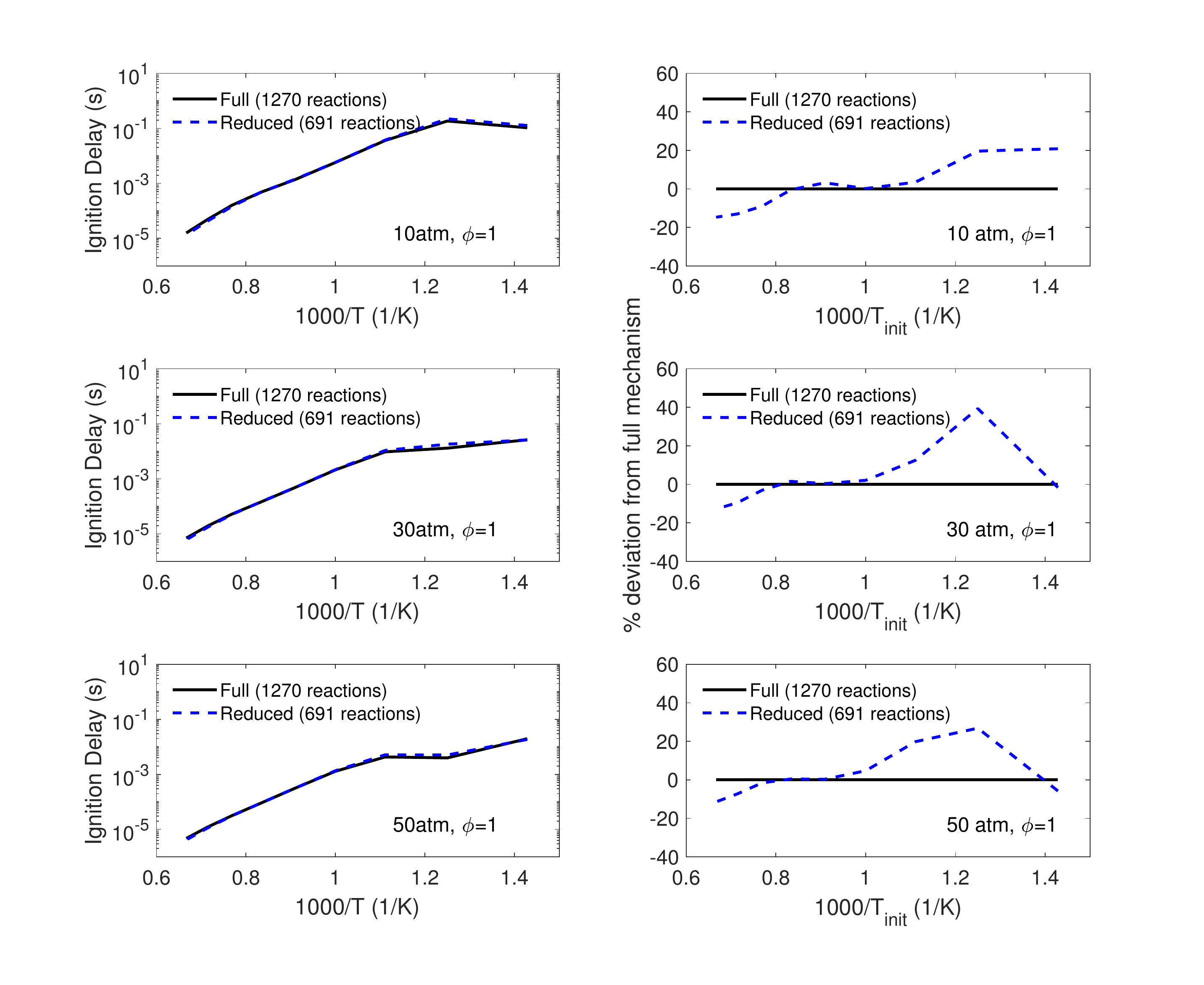}
\caption{Comparisons of computed ignition delay times using the full and reduced mechanisms for stoichiometric propane/air mixture at 10 atm, 30 atm, and 50 atm}
\label{fig:propane_ID}
\end{figure}

Using a similar approach as for the \ce{H2} mechanism study, we assembled a reduced mechanism for the propane chemistry.
To formulate a reduced mechanism that works in both low-/high\hyph temperature regime, we analyzed influential reactions identified from 216 conditions (700--1500~K, 1--50~atm, and equivalence ratios of 0.5--2) from homogeneous reactor simulations.
The performances of the reduced version of the propane mechanism, which consists of 111 species and 691 reactions, were again tested by comparing ignition delay times with the ones predicted by the parent mechanism.
The results, of which a representative sample is shown in Fig.~\ref{fig:propane_ID}, show that ignition delay times predicted by the reduced mechanism are within 40\% of those by the parent propane mechanism and that the NTC (negative temperature coefficient) behavior predicted by the parent mechanism is preserved.

\section{Conclusions}
\label{sec:conclusions}
In this work, we present a new method to analyze complex chemical reaction networks by employing a data-driven sparse-learning approach.
Using concentration and reaction rates at a given time, the proposed optimization approach is able to identify the most influential reactions in the network, that is the smallest subset of reactions that is able to approximate the concentration of the species to within a prescribed error tolerance using a convex relaxation method to solve the underlying optimization problem. The proposed approach has the following advantages: guaranteed concentration approximation accuracy over all time points; low computational cost; and versatility due to its applicability to general chemical reaction networks.  

We tested our approach on reaction networks generated by the combustion of two fuels, hydrogen and propane.
When applied to the \ce{H2} combustion, our method was able to identify the key reactions that mark the different phases of the ignition process.
Moreover, a reduced mechanism of the \ce{H2} oxidation is built by collecting the influential reactions at all times in a wide range of thermodynamic conditions (5--20~atm, 800--1100~K, $\phi=$~0.5--2), displayed a deviation from ignition delay time of less than 2\%, while using only half of the reactions.

The analysis of the \ce{C3H8} combustion mechanism showed that our method can identify the changes in low\hyph temperature pathways and capture the propane's characteristic two-stage ignition behavior, as well as the differences in relevant reactions between low\hyph ~and high\hyph temperature ignition conditions.
Similarly to hydrogen combustion, we built a reduced mechanism consisting of 111 species (reduced by 5.1\%) and 691 reactions (reduced by 45.6\%), by analyzing the combustions of 216 different systems.
The ignition delay times obtained with the reduced mechanism are within 40\% deviation of the original mechanism in a wide range of conditions (700--1500~K, 1--50~atm, $\phi=$~0.5-2).

This study showcases the potential of the proposed data-driven approach to analyze very complex reaction networks and to perform mechanism reduction in a computationally-efficient manner.

\section*{Acknowledgments}
This research has been funded by the US Army Research Office grants W911NF-15-1-0241 and W911NF-14-1-0359.
PE and AV thank the College of Engineering at the University of Michigan for partially supporting this work.

\bibliography{unified}

\bibliographystyle{model1-num-names}

\end{document}